\documentclass[letterpaper, 10 pt, conference]{ieeeconf}  % Comment this line out
\IEEEoverridecommandlockouts                              % This command is only
\usepackage{epsfig} % for postscript graphics files
\usepackage{amsmath,amssymb,amsfonts}
\usepackage{color}
\usepackage{balance}
\usepackage{booktabs,tabularx}
\usepackage{multirow}
\usepackage{mathtools}
\usepackage{pdfrender}
\usepackage{trfsigns}

\usepackage{centernot}
\usepackage{color}
\usepackage[acronym]{glossaries}
\usepackage[ruled]{algorithm2e}
\usepackage{graphicx,subcaption}
\usepackage{tabu}
\usepackage{cite}
\usepackage[official]{eurosym}

\graphicspath{{figures/}}

\newcommand{\R}{\mathbb{R}}

\newcommand{\N}{\mathbb{N}}

\newcommand{\mc}[1]{\mathcal{#1}}

\newcommand{\bs}{\boldsymbol}

\newtheorem{theorem}{Theorem}

\newtheorem{lemma}{Lemma}

\newtheorem{remark}{Remark}
\newtheorem{assumption}{Assumption}

\usepackage{xcolor,calc}

\makeglossaries
\makeglossaries
\newacronym{NOESC}{NOESC}{numerical optimization-based extremum seeking control}
\newacronym{ESC}{ESC}{extremum seeking control}

\title{\LARGE \bf
Extremum seeking control of a class of constrained nonlinear systems
}

\author{Shuai Yuan,	Filippo Fabiani and Simone Baldi% <-this % stops a space
\thanks{S.~Yuan is with School of Astronautics, Harbin Institute of Technology, China. F.~Fabiani is with the Department of Engineering Science, University of Oxford, United Kingdom. S.~Baldi is with School of Mathematics, Southeast University, China. (e-mail: {\tt \footnotesize{shuaiyuan@hit.edu.cn}}, {\tt \footnotesize{filippo.fabiani@eng.ox.ac.uk}}, {\tt \footnotesize{s.baldi@tudelft.nl}}). This work was partially supported by the Heilongjiang Postdoctoral Foundation under Grant No.LBH-Z19058 and the National Nature Science Foundation of China(NSFC) under Grant No.12002103, and through the Government’s modern industrial strategy by Innovate UK, part of UK Research and Innovation, under Project LEO (Ref. 104781).}%
}

\begin{document}

\maketitle
\thispagestyle{empty}
\pagestyle{empty}

%%%%%%%%%%%%%%%%%%%%%%%%%%%%%%%%%%%%%%%%%%%%%%%%%%%%%%%%%%%%%%%%%%%%%%%%
\begin{abstract}

%This paper studies the extremum seeking control (ESC) problem for a class of nonlinear systems subject to output constraints. 
%To steer the system to optimize a predefined, thought a-priori unknown, performance function, we propose a novel numerical optimization-based extremum seeking control (NOESC) design consisting of a constrained numerical optimization method and an inversion-based feedforward controller. In particular, a projected gradient descent algorithm is exploited to produce the state sequence to optimize the performance function, whereas a suitable boundary value problem accommodate the finite-time state transition between each two consecutive points of the state sequence. Compared to available NOESC methods, in addition to explicitly deal with output constraints, the contribution made by the proposed method is twofold: i) it can deal with \red{more general full-state dependence in the performance function}; ii) it can handle a more general class of nonlinear systems with possibly unstable internal dynamics \red{to the extremum of the performance function}. The effectiveness of the proposed ESC scheme is shown through extensive numerical simulations. 

This paper studies the extremum seeking control (ESC) problem for a class of constrained nonlinear systems. Specifically, we focus on a family of constraints allowing to reformulate the original nonlinear system in the so-called input-output normal form.
To steer the system to optimize a performance function without knowing its explicit form, we propose a novel numerical optimization-based extremum seeking control (NOESC) design consisting of a constrained numerical optimization method and an inversion-based feedforward controller. In particular, a projected gradient descent algorithm is exploited to produce the state sequence to optimize the performance function, whereas a suitable boundary value problem accommodates the finite-time state transition between each two consecutive points of the state sequence. Compared to available NOESC methods, the proposed approach i) can explicitly deal with output constraints; ii) the performance function can consider a direct dependence on the states of the internal dynamics; iii) the internal dynamics do not have to be necessarily stable. The effectiveness of the proposed ESC scheme is shown through extensive numerical simulations. 

\end{abstract}

%\begin{IEEEkeywords}
%	Extremum seeking control, constrained numerical optimization, finite-time state transition, output constraints
%\end{IEEEkeywords}

%%%%%%%%%%%%%%%%%%%%%%%%%%%%%%%%%%%%%%%%%%%%%%%%%%%%%%%%%%%%%%%%%%%%%%%%
\section{Introduction}
\Gls{ESC} has been consistently attracting research attention in the past two decades due to its ability to find an extremum of a performance function without knowing its explicit form. Therefore, it finds applications in many engineering fields, spanning from power optimization \cite{Ghaffari_TCST:2014}, robotic motion control \cite{Koropouli_CEP:2016}, source seeking \cite{Fu_TIE:2009}, etc. To date, various \gls{ESC} methods have been proposed by, e.g., incorporating dither signals and averaging methods \cite{Krstic_SLC:2000,Tan_Auto:2006}, sliding mode control \cite{Pan_IJC:2003, Yin_Auto:2014}, sample-data optimization \cite{khong2013unified, Nesic_Auto:2013}, and numerical optimization methods \cite{Zhang_TAC:2007, Ye_Auto:2016}. Since the first two techniques may result in undesired oscillating behaviors of the closed-loop system, the latter two methods naturally entail dither-free and chattering-free behaviors, and consequently are of theoretical and practical interest. 
%This work focuses on numerical optimization-based extremum seeking control. 

\Gls{NOESC} was initially proposed by the seminal work \cite{Zhang_TAC:2007}, which combined several numerical optimization schemes, such as the gradient descent or Newton's methods, and traditional state regulation techniques to find the extremum of the performance function. Some works followed that made use of this numerical optimization and state regulation framework. The \gls{NOESC}  proposed in \cite{Zhang_TAC:2007} is extended to state feedback linearizable systems with parametric uncertainties and input disturbances by integrating adaptive control techniques \cite{Zhang_Auto:2009}. Performance optimization problem of a special class of nonlinear systems with unmodeled dynamics and disturbances is considered, and a robust \gls{ESC} based on a conjugate gradient method and an extended state observer is proposed in \cite{Ye_Auto:2016} to find the minimum of the performance function. In \cite{Alick_IJCAS:2015}, a simplex gradient-based optimization is exploited to construct a derivative-free \gls{ESC} based on measurements of the performance function for linear systems. However, despite of some pioneering works in the field of \gls{NOESC} for linear systems or input-state linearizable systems, not much efforts have been made for input-output linearizable systems, mainly due to the difficulty to deal with the uncontrollable internal dynamics \cite{Zhang_TAC:2007}. In addition, the presence of system constraints is ubiquitous in practice, a fact that has not yet been considered in the literature body on \gls{NOESC}. 

In contrast with the aforementioned works, we focus on \gls{NOESC} for input-output linearizable systems with output constraints. Specifically, we consider all those constraints that enable the so-called input-output normal form for the original nonlinear system. The proposed \gls{NOESC} scheme consists of a constrained numerical optimization step and finite-time state transition, where the former gives rise to a state sequence towards the extremum of an a-priori unknown performance function, and the latter aims to design controllers driving the system to evolve along the state sequence. Based on the state sequence, the boundary conditions for the trajectory of the output and that of the internal dynamics are formulated. Successively, inspired by \cite{Graichen_Auto:2005}, we adopt an inversion-based feedforward control method with free parameters to fulfil finite-time state transitions, thus guaranteeing that the performance function is asymptotically minimized. In summary, the contribution made by the proposed method is threefold: i) it can explicitly deal with output constraints; ii) the performance function can include a direct dependence on the states of the internal dynamics, where iii) these latter are not assumed to be necessarily stable.

The paper is organized as follows: we present the considered \gls{ESC} problem for nonlinear systems and some preliminaries in \S II. Successively, in \S III we introduce one constrained numerical optimization method and the finite-time transition mechanism. Finally, our theoretical findings are corroborated in \S IV through a numerical example.

\smallskip
\emph{Notation: } $\mathbb{N}_{+}$,  $\mathbb{R}^{n}$, $\mathbb{R}^{n\times n}$ represent the sets of positive natural numbers, $n$-component real vectors, and $n$ by $n$ real matrices, respectively. The operator $\nabla \cdot$ is the gradient. Bold symbols represent vectors, while italics symbols denote scalars. The Lie derivative is defined as $ \mc{L}_{\bs{f}}^{i}h(\bs{x}) =  \partial (\mc{L}_{\bs{f}}^{i-1}h(\bs{x}))/\partial \bs{x} \cdot \bs{f}(\bs{x})$.

%%%%%%%%%%%%%%%%%%%%%%%%%%%%%%%%%%%%%%%%%%%%%%%%%%%%%%%%%%%%%%%%%%%%%%%%
\section{Problem formulation}
We consider the following nonlinear system in affine form with output constraints
\begin{equation} \label{eq:syst}
	\left\{
	\begin{aligned}
		\dot{\bs{x}}(t) = &\  \bs{f}(\bs{x}(t))+\bs{g}(\bs{x}(t))u(t), \; \bs{x}(t_{0}) = x_{0}, \\
		y(t) = &\ h(\bs{x}(t)), \; y \in \mc{Y},
	\end{aligned}
	\right.
\end{equation}
which has to be controlled to a minimizer of the following, \textit{a-priori} unknown performance function
\begin{equation} \label{eq:perfFunc}
	z(t) = J(\bs{x}(t)).
\end{equation}
Here, $\bs{x} \in \mathbb{R}^{n}$ denotes the state vector, $u \in \mathbb{R}$ the single input acting on the system, $y \in \mathbb{R}$ the system output with $\mc{Y}$ being a closed and convex set, and $z \in \mathbb{R}$ the performance output. We assume the map $h$ between the state and the output to be convex, and therefore, by making use of $\mc{Y}$, we can equivalently define a set of convex state constraints denoted by $\mc{X} \subseteq \R^n$ \cite{Boyd_book:2004}. Moreover, the mapping $\bs{f}:  \mc{X} \to \mathbb{R}^{n}$ and functions $h: \mc{X} \to \mathbb{R}$, $J:  \mc{X} \to \mathbb{R}$ are assumed to be continuously differentiable. Despite the fact that the explicit form of $J(\cdot)$ is not available, we assume to directly measure the value $J(\bs{x})$ through the variable $z$ in \eqref{eq:perfFunc}, whereas its gradient $\nabla J(\bs{x})$ can be estimated numerically \cite{Brekelmans_JOTA:2015}. 

%Thus, the overall goal of the \gls{ESC} scheme is to exploit the measurements of the output, the performance output, and the state to design a controller that drives the output-constrained system in \eqref{eq:syst} to a state that minimizes the unknown performance function $J(\cdot)$. In the spirit of traditional \gls{NOESC} schemes available in the literature, the proposed \gls{ESC} problem can be divided into two sub-problems: 
Thus, the goal of the \gls{ESC} scheme is to exploit the performance output and state measurements (the output can be reconstructed, accordingly) to design a controller that drives the output-constrained system in \eqref{eq:syst} to a state that minimizes the unknown function $J(\cdot)$. Following traditional \gls{NOESC} schemes available in the literature, the proposed \gls{ESC} problem can be divided into two sub-problems: 
\begin{enumerate}
	\item A constrained optimization problem, described by
	\begin{equation}\label{eq:min_prob}
		\begin{aligned}
			&\underset{\bs{x} \in \mc{X}}{\mathrm{min}} & & J(\bs{x}) \, .
		\end{aligned}
	\end{equation}
	By relying on available measurements, a minimizer of \eqref{eq:min_prob} is computed through an iterative scheme, which provides a state sequence $\left\{ \bs{x}_{k} \right\}_{k \in \mathbb{N}_{+}}$;
	
	\item A finite-time state transition problem, aiming to steer the system \eqref{eq:syst} from $\bs{x}_k$ to $\bs{x}_{k+1}$ within a finite time.
\end{enumerate}

%Before proposing the methods to solve the aforementioned two problems, 
We conclude by postulating the following standard assumption on the smoothness and convexity of the function $J(\cdot)$.

\smallskip
\begin{assumption} \label{assump1}
	The performance function $J$ is convex, and has $L$-Lipschitz continuous gradient, i.e.,
	\begin{equation*}
		\|\nabla J(\bs{x}) - \nabla J(\hat{\bs{x}}) \| \leq  L \|\bs{x}-\hat{\bs{x}}\|  \text{ for all } \bs{x}, \hat{\bs{x}} \in \mc{X}.
	\end{equation*}
%	 The unique minimizer $\bs{x}^{*} \in \mc{X}$ is such that $J(\bs{x}^{*}) < J(\bs{x})$, for all $\bs{x} \in \mc{X}$.
	\hfill$\square$
\end{assumption}
%\smallskip
%\begin{assumption}
%%	The function $J$ is $L$-smooth, i.e. 
%	The performance function $J$ has $L$-Lipschitz continuous gradient
%	\begin{equation*}
%		\|\nabla J(\bs{x}) - \nabla J(\hat{\bs{x}}) \| \leq  L \|\bs{x}-\hat{\bs{x}}\| \text{ for all } \bs{x}, \hat{\bs{x}} \in \mc{X}.
%	\end{equation*}
%	\hfill$\square$
%\end{assumption}

%%%%%%%%%%%%%%%%%%%%%%%%%%%%%%%%%%%%%%%%%%%%%%%%%%%%%%%%%%%%%%%%%%%%%%%%%
\section{Extremum seeking control}

In this section, we propose an \gls{ESC} design to solve the constrained optimization problem and to guarantee the finite-time state transition previously introduced.

\subsection{Constrained optimization}

%To iteratively minimize the performance function, constrained numerical optimization methods are applicable by making use of the measurements of $J(\bs{x})$ and its gradient estimate $\nabla J(\bs{x})$, which. Inspired by the numerical optimization-based extremum seeking introduced in \cite{Zhang_TAC:2007}, we propose to adopt constrained numerical optimization methods to find $J(\bs{x}^*)$. In this paper, the projected gradient optimization is exploited to solve the constrained numerical optimization problem.

Inspired by the \gls{NOESC} scheme introduced in \cite{Zhang_TAC:2007}, we adopt a projected gradient algorithm to iteratively compute a constrained state sequence that aims to seek $J(\bs{x}^*)$, where $\bs{x}^*$ denotes some minimizer of \eqref{eq:min_prob}, thus accommodating item 1) in \S II. Now, let us introduce the following key notion.

\smallskip
\begin{lemma} \textup{(Projection principle, \cite[Prop.~2.1.3]{bertsekas1997nonlinear})} \label{lem1}
	Let $\mc{X}$ be a closed and convex set. Then, $\bs{x}_{p}$ is the projected point of $\bs{z} \in \mathbb{R}^{n}$ onto $\mc{X}$ if and only if it satisfies
	\begin{equation*}
		(\bs{z}-\bs{x}_p)^\top(\bs{x}-\bs{x}_p) \leq 0, \text{ for all } \bs{x} \in \mc{X}.
	\end{equation*}
	\hfill$\square$
\end{lemma}
\smallskip

Therefore, we denote the projection of a point $\bs{z} \in \mathbb{R}^{n}$ onto a set $\mc{X}$ as
$
	\Pi_{\mc{X}}(\bs{z}) = \mathrm{argmin}_{\bs{x} \in \mc{X}} \, \tfrac{1}{2} \, \|\bs{x}-\bs{z} \|^2.
$
The main steps of the projected-gradient algorithm are summarized in Algorithm~\ref{alg:grad_desc}. Specifically, given a feasible starting point $\bs{x}_0 \in  \mc{X}$ and a (fixed) step-size $\varrho > 0$, the projected gradient descent algorithm produces a state sequence $\left\{ \bs{x}_k\right\}_{k=1}^{S_N}$, with $S_N \in \mathbb{N}_+$ defining a finite number of iterations, as follows
\begin{equation} \label{eq:PGD}
	\bs{x}_{k+1} = \Pi_{\mc{X}}\left( \bs{x}_{k} - \varrho \nabla J(\bs{x}_k) \right),
\end{equation}
where the step size $\varrho$ is selected to be $\frac{2}{L+2\epsilon}$ with $\epsilon > 0$, since the function $J$ is $L$-smooth, according to Assumption~\ref{assump1} \cite{Levitin_USSR:1966}.

\smallskip
\begin{remark}
	In determining the (fixed) step-size $\varrho$, for simplicity we assume to have available $L$, the Lipschitz constant of the gradient of the unknown performance function $J(\cdot)$. To remove the knowledge of $L$, we could adopt a backtracking line search, or variable step-sizes, i.e., $\varrho = \varrho_k$, to guarantee the convergence of the state sequence \cite{bertsekas1997nonlinear,nocedal2006numerical}.
	\hfill$\square$
\end{remark}

\begin{algorithm}[!t]
	\caption{Projected gradient descent method}\label{alg:grad_desc}
	\DontPrintSemicolon
	\SetArgSty{}
	\SetKwFor{ForAll}{for all}{do}{end forall}
	\smallskip
	\textbf{Initialization: } Choose a starting point $\bs{x}_{0}$, set $k \coloneqq 0$, $\varepsilon_{0} > 0$, and $\varrho \coloneqq 2/(L+2\epsilon)$, for some $\epsilon > 0$\\
	\smallskip
	\While{$\|\nabla J(\bs{x}_k))\| \geq \varepsilon_{0}$}{
		\begin{itemize}\setlength{\itemindent}{-.5cm}
			\smallskip
			\item[] $\hat{\bs{x}}_{k+1} \longleftarrow \bs{x}_{k} - \varrho \nabla J(\bs{x}_k) $
			\smallskip
			\item[] $\bs{x}_{k+1} \longleftarrow \Pi_{\mc{X}} (\hat{\bs{x}}_{k+1}) $
			\smallskip
			\item[] $k \longleftarrow k+1$
		\end{itemize}
	}
\end{algorithm}

%\begin{algorithm}
%	\caption{Projected Gradient Descent}
%	\begin{algorithmic}
%		\textbf{Initial:} Given $\bs{x}_{0}$, $\varepsilon_{0}$, $k=0$ \newline
%		\WHILE{$\|\nabla J(\bs{x}_k))\| \geq \varepsilon_{0}$}
%		\STATE $\varrho_k \longleftarrow \frac{2}{L+2\epsilon}$
%		\STATE $\hat{\bs{x}}_{k+1} \longleftarrow \bs{x}_{k} - \varrho_k \nabla J(\bs{x}_k) $
%		\STATE $\bs{x}_{k+1} \longleftarrow \Pi_{\mc{X}} (\hat{\bs{x}}_{k+1}) $
%		\STATE $k \longleftarrow k+1$
%		\ENDWHILE
%	\end{algorithmic}
%\end{algorithm}

\subsection{Finite-time state transition}
In this part, the state sequence $\left\{ \bs{x}_k\right\}_{k\in \mc{S}}$ produced by the numerical optimization algorithm is exploited to formulate the boundary conditions for system \eqref{eq:syst}. Thus, by introducing the discrete time instant $t_{k}$, $k\in \mc{S} \coloneqq \left\{1,2,\dots,S_N\right\}$, and
%$
%	t_{1},t_{2},\dots, \dots, \, k\in \mc{S} := \left\{1,2,\dots,S_N\right\} 
%$
$\Delta_k \coloneqq t_{k+1}-t_{k}$ as a tuning parameter, we define the two-point boundary conditions for the state variable $\bs{x}$ as follows
\begin{equation} \label{eq:stateBC} 
	\bs{x}(t_{k}) = \bs{x}_{k}, \; \bs{x}(t_{k+1}) = \bs{x}_{k+1}, \; k\in \mc{S},
\end{equation}
where $\bs{x}_{k+1}$ is obtained by making use of the constrained numerical optimization as in Algorithm 1. 
From a mathematical point of view, the $n$ ordinary differential equations of \eqref{eq:syst} and the $2n$ boundary conditions \eqref{eq:stateBC}  form a two-point boundary value problem for the states $\bs{x}(t)$, $t \in [t_{k},t_{k+1}]$.
In this context, the finite-time transition problem aims to design a state-feedback controller that ensures the transition between every two consecutive points $\bs{x}_{k}$ and $\bs{x}_{k+1}$ of the system within a finite time interval $t \in [t_{k},t_{k+1}]$. 

Now, by leveraging the system model in \eqref{eq:syst}, we can translate the boundary conditions on the state variable in \eqref{eq:stateBC} into some boundary conditions on the system output, i.e.,
\begin{equation} \label{eq:outputBC}
	y(t_{k}) = y_k = h(\bs{x}_{k}), \; y(t_{k+1}) = y_{k+1} = h(\bs{x}_{k+1}), \; k \in \mc{S}.
\end{equation}

In this work, we propose to extend the inversion-based feedforward control method proposed in \cite{Graichen_Auto:2005} based on the input-output normal form, so as to achieve finite-time transitions between the state sequence $\left\{\bs{x}_{k} \right\}_{k\in \mc{S}}$. To this end, let us assume that the system in \eqref{eq:syst} has a well-defined relative degree $r$, with $0 < r \leq n$ \footnote{The relative degree $r$ means how often the output $y$ has to be differentiated when the input $u$ appears explicitly -- see, e.g., \cite{isidori_Book:1995}.}. Then, by introducing normal coordinates $y, \, \dot{y}, \,\dots, \,y^{(r-1)}$, the system \eqref{eq:syst} can be transformed into the following input-output normal form \cite{Nijmeijer_Book:1990}
\begin{equation} \label{eq:IOsys}
	y^{(r)} = \alpha \left(y, \dot{y},\dots,y^{(r-1)}, \bs{\eta}, u \right)
\end{equation}
with the internal dynamics
\begin{equation} \label{eq:InterSyst}
	\dot{\bs{\eta}} = \bs{\beta} \left(y, \dot{y},\dots,y^{(r-1)}, \bs{\eta}, u \right)
\end{equation}
where $\bs{\eta}:=\left[\eta_{r+1},\dots,\eta_{n} \right]^\top \in \mathbb{R}^{n-r}$ denotes the state of the internal dynamics, while the mappings $\alpha:\mathbb{R}^{n}\times \mathbb{R} \to \mathbb{R}$ and $\bs{\beta}:\mathbb{R}^{n} \times \mathbb{R} \to \mathbb{R}^{n}$ depend on the original system \eqref{eq:syst}.
Thus, the new coordinates $y, \, \dot{y}, \,\dots, \,y^{(r-1)}$ and internal dynamics $\bs{\eta}$ are combined to complete the diffeomorphism \cite{Slotine_Book:1990}
\begin{equation}\label{eq:diffeo}
	\begin{aligned}
		&\ \left[y,\dot{y},\dots,y^{(r-1)}, \bs{\eta}^\top \right]^\top \\
		&\ \quad \quad \quad  = \left[ \phi_{1}(\bs{x}), \dots, \phi_{r+1}(\bs{x}), \dots,\phi_{n}(\bs{x}) \right]^\top = \bs{\phi}(\bs{x}) 
	\end{aligned}
\end{equation}
where the output derivatives are computed by
\begin{equation*}
	\left\{
	\begin{aligned}
		&\ y^{(i)} =  \mc{L}_{\bs{f}}^{i}h(\bs{x}) = \phi_{i+1}(\bs{x}), \; i=0,\dots,r-1, \\
		&\ \frac{\partial \phi_{j}}{\partial \bs{x}} \bs{g} = 0, \; j = r+1,\dots,n.   
	\end{aligned}
	\right.
\end{equation*}

\smallskip
\begin{remark}
	The coordinate $\bs{\eta}$ shall be selected such that the Jacobian of the diffeomorphism $\bs{\phi}(\bs{x})$ is invertible, thus guaranteeing the existence of its inverse $\bs{\phi}^{-1}$. For instance, a well-known method to find such $\bs{\eta}$ is to solve the partial differential equations $\nabla \eta_{j}^\top \bs{g} = 0$, $j = r+1, \dots, n$ \cite{Slotine_Book:1990}.
	\hfill$\square$
\end{remark}
\smallskip

Therefore, in view of the diffeomorphism in \eqref{eq:diffeo}, the boundary conditions for the internal dynamics read as
\begin{equation}\label{eq:etaBC}
	\left\{
	\begin{aligned}
		&\ \bs{\eta}(t_{k}) = \bs{\eta}_{k} = \left[\phi_{r+1}(\bs{x}_{k}),\dots, \phi_{n}(\bs{x}_{k})\right]^\top, \\
		&\ \bs{\eta}(t_{k+1}) = \bs{\eta}_{k+1} = \left[\phi_{r+1}(\bs{x}_{k+1}),\dots, \phi_{n}(\bs{x}_{k+1})\right]^\top.
	\end{aligned}
	\right.
\end{equation}
According to the input-output normal form in \eqref{eq:IOsys}--\eqref{eq:InterSyst}, the inversion-based feedforward controller is designed by reverting the relation in \eqref{eq:IOsys}, and therefore it turns out to be
\begin{equation} \label{eq:InvControl0}
	u = \alpha^{-1}\left(y,\dots,y^{(r)},\bs{\eta} \right),
\end{equation}
thus revealing an explicit dependence on the output trajectory $y$ and the state $\bs{\eta}$ of the internal dynamics.  Thus, to obtain explicitly the controller \eqref{eq:InvControl0}, let us first select an a-priori reference output trajectory for the system output $y$, say $y^*$, that is at least $r$ times differentiable and that satisfies the two-point boundary conditions in \eqref{eq:outputBC}, i.e., for all $k \in \mc{S}$,
\begin{equation} \label{eq:outputBC1}
	y^*(t_{k}) = y_{k}, \; y^{*}(t_{k+1}) = y_{k+1}
\end{equation}
which recasts the input-output normal form into:
\begin{equation}\label{eq:OIsyst}
	\left\{
	\begin{aligned}
		y^{*(r)} = \alpha \left(y^*,\dot{y}^*,\dots,y^{*(r-1)},\bs{\eta}^*,u^* \right), \\
		\dot{\bs{\eta}}^* = \bs{\beta} \left( y^*,\dot{y}^*,\dots,y^{*(r-1)},\bs{\eta}^*,u^* \right).
	\end{aligned}
	\right.
\end{equation}
Here, $u^*$ and $\eta^*$ represent the input and the state of the internal dynamics corresponding to the pre-selected output trajectory $y^*$, respectively.
Then, according to \cite{Graichen_Auto:2005}, the inversion-based feedforward control design follows readily from the input-output form the relation in \eqref{eq:OIsyst} as
\begin{equation} \label{eq:InvContrl}
	u^* = \alpha^{-1}\left(y^*,\dots,y^{*(r)},\bs{\eta}^* \right),
\end{equation}
and $\bs{\eta}^*$ should satisfy
\begin{equation} \label{eq:IntDynamics}
	\left\{
	\begin{aligned}
		&\dot{\bs{\eta}}^* =\bs{\beta} \left(y^*, \dots,y^{(r-1)^*}, \bs{\eta}, u^* \right) = \bs{\lambda} \left( y^*,\dots,y^{*(r-1)},\bs{\eta}^* \right) \\
		&\ \text{subject to } \bs{\eta}^*(t_{k}) = \bs{\eta}_{k}, \; \bs{\eta}^*(t_{k+1}) = \bs{\eta}_{k+1}.
	\end{aligned}
	\right.
\end{equation} 
In this case, the boundary conditions $\bs{\eta}_{k}$, for all $k\in \mc{S}$, are defined as in \eqref{eq:etaBC}, and  the function $\bs{\lambda}$ is obtained by replacing the controller \eqref{eq:InvContrl} into the internal dynamics described in \eqref{eq:OIsyst}. We note that the reference output trajectory $y^*$ and its derivatives are deemed as input of \eqref{eq:IntDynamics}. With this regard, the finite-time state transition problem reduces to the boundary value problem of \eqref{eq:IntDynamics} subject to the boundary constraints \eqref{eq:etaBC}, with a properly pre-selected output trajectory $y^*$. In other words, the inversion-based feedfordward controller fulfilling the finite-time state transition can be designed by exploiting the trajectory of the internal dynamics via \eqref{eq:InvContrl}.  

In what follows, we provide a systematic way to design the reference output trajectory $y^*$. First of all, we note that, in view of the fact that the function $h$ returns a scalar output, the constraint set $\mc{Y}$ amounts to box constraints, i.e., $y^{*} \in \left[y_{\min}, y_{\max} \right]$, where $y_{\min} < y_{\max}\in \mathbb{R}$ denote a lower and an upper bound of the output, respectively. Additionally, $y^*$ shall also satisfy the two-point boundary conditions in \eqref{eq:outputBC1}, for all $k \in \mc{S}$.
For these reasons, we first aim to design an unconstrained trajectory $\zeta$, and then map it into a constrained reference output trajectory. Specifically, to fulfill the output constraints, we adopt a saturation function, and therefore the reference output trajectory can be defined as follows
\begin{equation} \label{eq:refY}
	y^* = \varphi \left(\zeta,y_{\min}^{\mathrm{s}},y_{\max}^{\mathrm{s}} \right),
\end{equation}
where $(y_{\min}^{\mathrm{s}}, y_{\max}^{\mathrm{s}})$ are the asymptotic bounds of the saturation function as $\|\zeta\| \to \infty$, formally defined as $\left[y_{\min}-\delta_{y}, y_{\max}+\delta_{y} \right]$ with $\delta_{y} > 0$. In addition, inspired by \cite{Graichen_TAC:2008}, $\zeta$ is a virtual and unconstrained output trajectory and the saturation function $\varphi$ is defined as a sigmoid function 
\begin{equation*}
	\begin{aligned}
		&\ \varphi(\zeta,y_{\min}^{\mathrm{s}},y_{\max}^{\mathrm{s}}) = y_{\max}^{\mathrm{s}} - \frac{y_{\max}^{\mathrm{s}} - y_{\min}^{\mathrm{s}}}{1+\exp\left( \rho \zeta \right)},\; \rho =\frac{4}{y_{\max}^{\mathrm{s}} - y_{\min}^{\mathrm{s}}}.
	\end{aligned}
\end{equation*}

Then, the boundary conditions for $\zeta$ can be written as
\begin{equation*}
	\begin{aligned}
		\zeta(t_{\kappa}) = &\ \zeta_\kappa = \varphi^{-1} \left(y^*, y_{\min}^{\mathrm{s}}, y_{\max}^{\mathrm{s}} \right) \\
		= &\ \frac{1}{\rho} \left[   \log( y^{*}_{\kappa} - y_{\min}^{\mathrm{s}}) - \log(y_{\max}^{\mathrm{s}} - y^{*}_{\kappa})\right],
	\end{aligned}
\end{equation*}
%where the boundary conditions should satisfy the 
which, on the other hand, are subject to the output constraints $y_{\max}^{\mathrm{s}} -y^{*}_{\kappa} > 0$ and $y^{*}_{\kappa} - y_{\min}^{\mathrm{s}} < 0$, for $\kappa = k,k+1$.

According to \eqref{eq:IntDynamics}, the boundary value problem of the internal dynamics is overdetermined by $2(n-r)$ boundary conditions for $(n-r)$ first-order ODEs. Therefore, for the solvability of the boundary value problem of the internal dynamics $\bs{\eta}$, the following virtual output trajectory defined by a two-order ansatz function incorporating $n-r$ free parameters $\bs{p} = \left(p_{1},\dots,p_{n-r} \right)$ is proposed
\begin{equation} \label{eq:zetaFunc}
	\begin{aligned}
		\zeta(t,\bs{p}) = \zeta_k + &\ \sum_{i=r-1}^{2r-1} a_{i}(\gamma_i p_i,\zeta_{k},\zeta_{k+1}) \left( \frac{t-t_{k}}{\Delta_k}\right) \\
		&\ \quad \quad \quad \quad \quad  + \sum_{i=1}^{n-r} \gamma_i p_{i} \left(\frac{t-t_{k}}{\Delta_k} \right)^2,
	\end{aligned}
\end{equation}
for $t \in [t_k, \, t_{k+1}]$, and $\gamma_{i} > 0$ is a tuning parameter. Moreover, the linear functions $a_i$ are chosen so that the conditions $\zeta(t_{k},\bs{p}) = \zeta_{k}$, and $\zeta(t_{k+1},\bs{p}) = \zeta_{k+1}$ are met.
%\begin{equation*}
%	\zeta(t_{k},\bs{p}) = \zeta_{k}, \quad \zeta(t_{k+1},\bs{p}) = \zeta_{k+1}.
%\end{equation*}

\smallskip
\begin{remark}
	It is known that two boundary points can be satisfied with a one-order function that defines a straight line. However, a straight-line function is not able to accommodate for $n-r$ free parameters.  This is essentially the reason that motivates us to propose the second-order ansatz function in \eqref{eq:zetaFunc} to solve the two-boundary value problem in \eqref{eq:IntDynamics}.  Note that the ansatz function \eqref{eq:zetaFunc} significantly differs with the one developed in \cite{Graichen_TAC:2008}: a tuning parameter $\gamma_i$ is introduced to adjust the shape of the virtual output trajectory.
	\hfill$\square$
\end{remark}
\smallskip

In summary, the finite-time state transition reduces to  the two-point boundary value problem for the internal dynamics $\bs{\eta}$ with the design parameter $\bs{p}$
\begin{equation} \label{eq:etaStar}
	\dot{\bs{\eta}}^* = \bs{\psi} \left( \varphi(\bs{p},y_{\min}^{\mathrm{s}},y_{\max}^{\mathrm{s}}),\dots,\varphi^{(r-1)}(\bs{p},y_{\min}^{\mathrm{s}},y_{\max}^{\mathrm{s}}),\bs{\eta}^* \right) 
\end{equation}
with boundary conditions
$
%\begin{equation*}
	\bs{\eta}^*(t_{k}) = \bs{\eta}_{k},$ $\bs{\eta}^*(t_{k+1}) = \bs{\eta}_{k+1},
%\end{equation*}
$
where $\bs{\psi}$ follows by substituting \eqref{eq:zetaFunc} and \eqref{eq:refY} into \eqref{eq:IntDynamics}.  The proposed \gls{NOESC} scheme is summarized in Algorithm~\ref{alg:noesc}, while its convergence properties are established next.

\begin{algorithm}[!t]
	\caption{\gls{NOESC} for output constrained, nonlinear systems}\label{alg:noesc}
	\DontPrintSemicolon
	\SetArgSty{}
	\SetKwFor{ForAll}{for all}{do}{end forall}
	\smallskip
	\textbf{Initialization: } Choose a starting point $\bs{x}_{0}$, set $k \coloneqq 0$, $\varepsilon_{0} > 0$, $t_{0}$\\
	\smallskip
	\While{$\|\nabla J(\bs{x}(t_{k+1}))\| \geq \varepsilon_{0}$}{
		\begin{itemize}\setlength{\itemindent}{-.5cm}
			\smallskip
			\item[] Obtain $\bs{x}_{k+1}$ via Algorithm~\ref{alg:grad_desc};
			\smallskip
			\item[] Define $y_{k}$, $y_{k+1}$, $\bs{\eta}_{k}$, and $\bs{\eta}_{k+1}$ via \eqref{eq:outputBC} and \eqref{eq:etaBC}; 
			\smallskip
			\item[] Design $\zeta$ via \eqref{eq:zetaFunc};
			\smallskip
			\item[] Obtain $\bs{\eta}^*$ by solving \eqref{eq:IntDynamics};
			\smallskip
			\item[] $k \longleftarrow k+1$;
		\end{itemize}
	}
\end{algorithm}

%\begin{algorithm}
%	\caption{NOESC for nonlinear systems with output constraints}
%	\begin{algorithmic}
%		\textbf{Initial:} Given $\bs{x}_{0}$, $\epsilon_{0}$, set $t_{0}$ and $k=0$ \newline
%		\WHILE{$\|\nabla J(\bs{x}(t_{k+1}))\| \geq \epsilon_{0}$}
%		\STATE Obtain $\bs{x}_{k+1}$ via \textbf{Algorithm 1};
%		\STATE Define $y_{k}$, $y_{k+1}$, $\bs{\eta}_{k}$, and $\bs{\eta}_{k+1}$ via \eqref{eq:outputBC} and \eqref{eq:etaBC}; 
%		\STATE Design $\zeta$ via \eqref{eq:zetaFunc};
%		\STATE Obtain $\bs{\eta}^*$ by solving \eqref{eq:IntDynamics};
%		\STATE Obtain $u^*(t)$ as $t\in [t_{k},t_{k+1}]$;
%		\STATE $k \longleftarrow k+1$;
%		\ENDWHILE
%	\end{algorithmic}
%\end{algorithm}

\begin{figure}[t!]
	\centering
	\includegraphics[width=1\columnwidth]{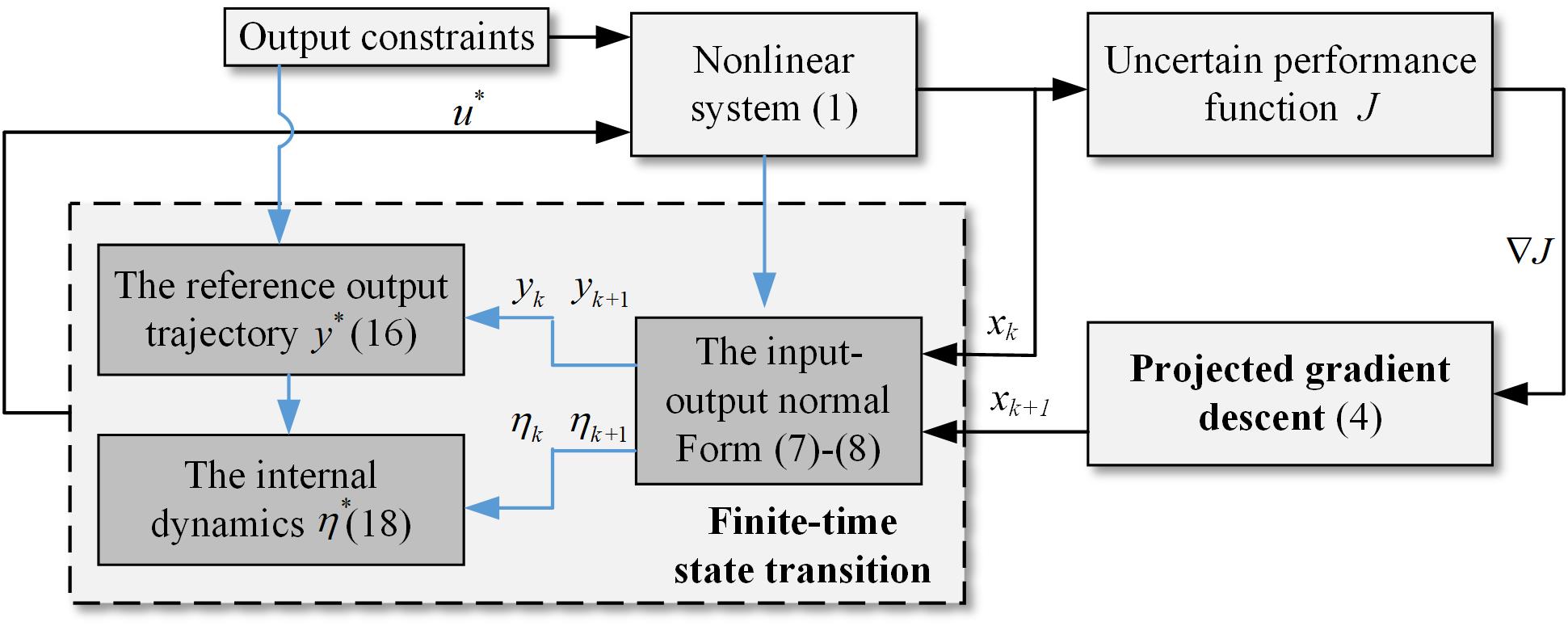}
	\caption{The proposed NOESC mechanism.}
	\label{fig:controlScheme}
\end{figure}

\smallskip
\begin{theorem}
%	Consider the nonlinear system \eqref{eq:syst} with the performance function \eqref{eq:perfFunc}. Under Assumption~1 and 2, suppose that the step size as in \textbf{Algorithm~1} is chosen and that there exists a trajectory $\bs{\eta}^*$ satisfying \eqref{eq:etaStar}. Then 
%	the state sequence $\{\bs{x}_k\}_{k \in \N}$ generated by the \gls{NOESC} scheme described in Algorithm~\ref{alg:noesc}
%	the state will 
%	globally asymptotically converges to the unique minimizer of the performance function in \eqref{eq:perfFunc}.
	Let Assumption~\ref{assump1} hold true. If, in addition, a trajectory $\bs{\eta}^*$ satisfying \eqref{eq:etaStar} exists, then the \gls{NOESC} scheme, described in Algorithm~\ref{alg:noesc}, generates a state sequence $\{\bs{x}_k\}_{k \in \N}$ for \eqref{eq:syst} that globally asymptotically converges to a minimizer of the performance function in \eqref{eq:perfFunc}.
	\hfill$\square$
\end{theorem}
\smallskip
\begin{proof}
	Let $\varrho$ be chosen as in Algorithm~\ref{alg:grad_desc}, i.e., $\varrho = 2/(L+2\epsilon)$, for some $\epsilon > 0$.
	In view of Assumption~\ref{assump1}, $\nabla J(\bs{x})$ is Lipschitz continuous with constant $L$, and hence we have
	\begin{equation*}
		\begin{aligned}
			J(\bs{x}_{k+1}) &\ - J(\bs{x}_k) \\
			\leq &\ \left\langle \bs{x}_{k+1}-\bs{x}_k, \nabla J(\bs{x}_k) \right\rangle + \frac{L}{2} \|\bs{x}_{k+1}-\bs{x}_k\|^2 \\
			= &\ -\frac{1}{\varrho}\left\langle \bs{x}_{k}-\varrho \nabla J(\bs{x}_k)-\bs{x}_{k+1}, \bs{x}_{k+1}-\bs{x}_k \right\rangle \\
			&\ - \frac{1}{\varrho} \|\bs{x}_{k+1}-\bs{x}_k\|^2 + \frac{L}{2}\|\bs{x}_{k+1}-\bs{x}_k\|^2.
		\end{aligned}
	\end{equation*}
	By leveraging the projection principle in Lemma~\ref{lem1}, we obtain
	\begin{equation*}
		\left\langle \bs{x}_{k}-\varrho \nabla J(\bs{x}_k)-\bs{x}_{k+1}, \bs{x}_{k+1}-\bs{x}_k \right\rangle \geq 0
	\end{equation*}
	which results in
	\begin{equation} \label{eq:JDiff}
		\begin{aligned}
			J(\bs{x}_{k+1}) - J(\bs{x}_k) \leq &\  \left( - \frac{1}{\varrho} + \frac{L}{2} \right)  \|\bs{x}_{k+1}-\bs{x}_k\|^2 \\
			\leq &\ - \epsilon \|\bs{x}_{k+1}-\bs{x}_k\|^2,
		\end{aligned}
	\end{equation}
%	for a fixed step $\varrho_k = \varrho = 2/(L+2\epsilon)$, for some $\epsilon > 0$. 
	The relation in \eqref{eq:JDiff} directly implies that the performance function $J(\bs{x}_k)$ is monotonically decreasing and, in addition, that $J(\bs{x}_{k+1})-J(\bs{x}_k) \to 0$. Equivalently, we have that $\bs{x}_{k+1} - \bs{x}_{k} \to 0$ for $k \to \infty$, as $\|\bs{x}_{k+1}-\bs{x}_k\|^2 \leq (J(\bs{x}_{k}) - J(\bs{x}_{k+1}))/\epsilon$. Therefore, according to Assumption~\ref{assump1}, we have 
	\begin{equation*}
		\begin{aligned}
			&\ J(\bs{x}_k) -  J(\bs{x}^*) \leq \left\langle \bs{x}_{k} - \bs{x}^*, \nabla J(\bs{x}_k) \right\rangle  \\
			= &\ \left\langle \nabla J(\bs{x}_k), \bs{x}_k - \bs{x}_{k+1} \right\rangle  - \frac{1}{\varrho}\left\langle \bs{x}_{k}- \bs{x}_{k+1}, \bs{x}^*-\bs{x}_{k+1} \right\rangle \\
			&\  + \frac{1}{\varrho}\left\langle \bs{x}_{k}-\varrho \nabla J(\bs{x}_k)-\bs{x}_{k+1}, \bs{x}^*-\bs{x}_{k+1} \right\rangle \\
			\leq &\ \left( \|\nabla J(\bs{x}_k) \| + \|\bs{x}^* - \bs{x}_{k+1} \| \right) \|\bs{x}_{k+1} - \bs{x}_k\|.
		\end{aligned}
	\end{equation*}
	Thus, in view of the fact that $\|\bs{x}_{k+1} - \bs{x}_k\| \to 0$ as $k \to \infty$, we obtain that $J(\bs{x}_k) \to J(\bs{x}^*)$ as $k \to \infty$. 
	
	Finally, since we assume a trajectory $\bs{\eta}^*$ satysfying the boundary value problem \eqref{eq:etaStar} exists, the inversion-based feedforward  controller \eqref{eq:InvContrl} guarantees that 
	\begin{equation} \label{eq:bound1}
		\left\{
		\begin{aligned}
			&\ \bs{\eta}^*(t_{k})  = \left[\phi_{r+1}(\bs{x}_{k}),\dots, \phi_{n}(\bs{x}_{k})\right]^\top, \\
			&\ \bs{\eta}^*(t_{k+1})  = \left[\phi_{r+1}(\bs{x}_{k+1}),\dots, \phi_{n}(\bs{x}_{k+1})\right]^\top.
		\end{aligned}
		\right.
	\end{equation}
	Furthermore, the reference output trajectory $y^*$ is defined to satisfy the boundary condition for all $k \in \mc{S}$, i.e., 
	\begin{equation}\label{eq:bound2}
		y^*(t_{k}) = h(\bs{x}_k), \; y^*(t_{k+1}) = h(\bs{x}_{k+1}).
	\end{equation}
	Applying the inverse of the diffeomorphism $\phi(\bs{x})$ to the boundary conditions \eqref{eq:bound1}--\eqref{eq:bound2} leads to a finite-time state transition via the controller \eqref{eq:InvContrl}, which entails that the state of the system \eqref{eq:syst} asymptotically converges to a point that minimizes the performance function \eqref{eq:perfFunc}.
\end{proof}
\smallskip
\begin{remark}
	The solution of two–point boundary value problem with free parameters as given in \eqref{eq:etaStar} for the state $\bs{\eta}$ and the parameter set $\bs{p}$ can be efficiently obtained by various numerical methods (e.g., the function \textit{bvp4c} in MATLAB). Thus, the state and reference output trajectories, $\bs{\eta}^*$ and $y^*$, can be obtained based on the saturated function $\varphi$, which results in the desired controller $u^*$ that steers the system from $\bs{x}_k$ to $\bs{x}_{k+1}$ for $t \in \left[t_{k},t_{k+1} \right]$, as shown in Fig.~\ref{fig:controlScheme}.
	\hfill$\square$
\end{remark}

%%%%%%%%%%%%%%%%%%%%%%%%%%%%%%%%%%%%%%%%%%%%%%%%%%%%%%%%%%%%%%%%%%%%%%%%%
\section{Numerical example}
We run our numerical simulations by making use of a nonlinear dynamics, re-adapted from \cite{Graichen_Auto:2005}, i.e., 
\begin{equation}\label{eq:systExp}
	\left\{
	\begin{aligned}
		\dot{x}_{1} = &\ -x_{2}^3 + u, \\
		\dot{x}_{2} = &\ \rho \left( 2x_{1}^2 - 2x_{2} \right), \\
		y =&\ x_1,
	\end{aligned}
	\right.
\end{equation}
with the output constraint $-1.5 \leq y \leq 1.5$, and the parameter $\rho$ determines the stability of the internal dynamics, as clarified next.  The performance function is taken as
\begin{equation*}
	J(x_{1},x_{2}) = 100\left(x_{2} - x_{1}^2 \right)^2 + \left(1 - x_{1}  \right)^2,
\end{equation*}
which admits a unique minimum $J_{\min} = 0$ when $x_{1}=1$ and $x_{2}=1$.
By adopting the change of coordinate $y=x_{1}$ and $\eta =x_{2}$, the relative degree is $r=1$, and the input-output form with the internal dynamics of the original system \eqref{eq:systExp} reads as (the bold notation is dropped since $\eta$ is a scalar)
\begin{equation}\label{eq:IOsystExp}
	\left\{
	\begin{aligned}
		\dot{y} = &\ -\eta^3 + u =: \alpha(\eta,u), \\
		\dot{\eta} = &\ \rho \left( 2y^2 - 2\eta \right) =: \beta (y,\eta).
	\end{aligned}
	\right.
\end{equation}

Note that, for $\rho > 0$, the internal dynamics are stable, while for $\rho < 0$, are unstable. Next, different values of $\rho$ are selected to show the effectiveness of the proposed \gls{NOESC} scheme, regardless of the stability of the internal dynamics. 

%\vspace{0.1cm}
%\textit{A. Control performance with different selections of $\rho$} 
%\vspace{0.1cm}

\subsection{Control performance with different selections of $\rho$}

At every iteration of Algorithm~\ref{alg:noesc}, we have to face with the following two point boundary conditions on $y$ and $\eta$
\begin{equation*}
	\left\{
	\begin{aligned}
		y(t_{k}) = &\ x_{1,k}, \; y(t_{k+1}) = x_{1,k+1}, \\
		\eta(t_{k}) = &\ x_{2,k}, \; \eta(t_{k+1}) = x_{2,k+1}.
	\end{aligned}
	\right.
\end{equation*}

By selecting the asymptotic bound $(-2,2)$ for the function $\varphi$, we design the reference output trajectory $y^*$ as
\begin{equation} \label{eq:refOut}
	y^*(t) = 2 - \frac{4}{1+\exp(4\zeta(t,p))},
\end{equation}
with the associated virtual output trajectory
\begin{equation*}
	\zeta(t,p) = \zeta_{k} + \left(\zeta_{k+1} - \zeta_{k} - \gamma p \right) \left(\frac{t-t_k}{\Delta_k}\right) + \gamma p\left(\frac{t-t_k}{\Delta_k} \right)^2.
\end{equation*}
Here, the linear function $a = \zeta_{k+1} - \zeta_{k} - \gamma p$ is selected such that $\zeta(t_k,p) = \zeta_k$ and $\zeta(t_{k+1},p) = \zeta_{k+1}$.
In accordance, the boundary conditions for the virtual output trajectory are 
%defined by
\begin{equation*}
	\left\{
	\begin{aligned}
		\zeta_{k} = &\ \frac{1}{4}\left[ \log \left(y_{k} +2 \right) - \log \left( 2- y_{k}  \right) \right] \\
		\zeta_{k+1} =&\ \frac{1}{4}\left[ \log \left(y_{k+1} +2 \right) - \log \left( 2 - y_{k+1}\right) \right].
	\end{aligned}
	\right.
\end{equation*}
With $\rho = 1$ the internal dynamics are stable, and therefore with the reference output trajectory $y^*$ in \eqref{eq:refOut}, the updated internal dynamics and the inversion-based controller reads as
\begin{equation}
	\left\{
	\begin{aligned}
		\dot{\eta}^* = &\ 2y^* - 2 \eta^{*3}, ~\eta^{*}(t_{k}) = x_{2,k},~ \eta^{*}(t_{k+1})=x_{2,k+1} \\
		u^* = &\ \dot{y}^*  + \eta^{*3} =: \alpha^{-1} (\dot{y}^*,\eta^*).
	\end{aligned}
	\right.
\end{equation}

\begin{figure}
	\centering
	\includegraphics[width=0.5\textwidth]{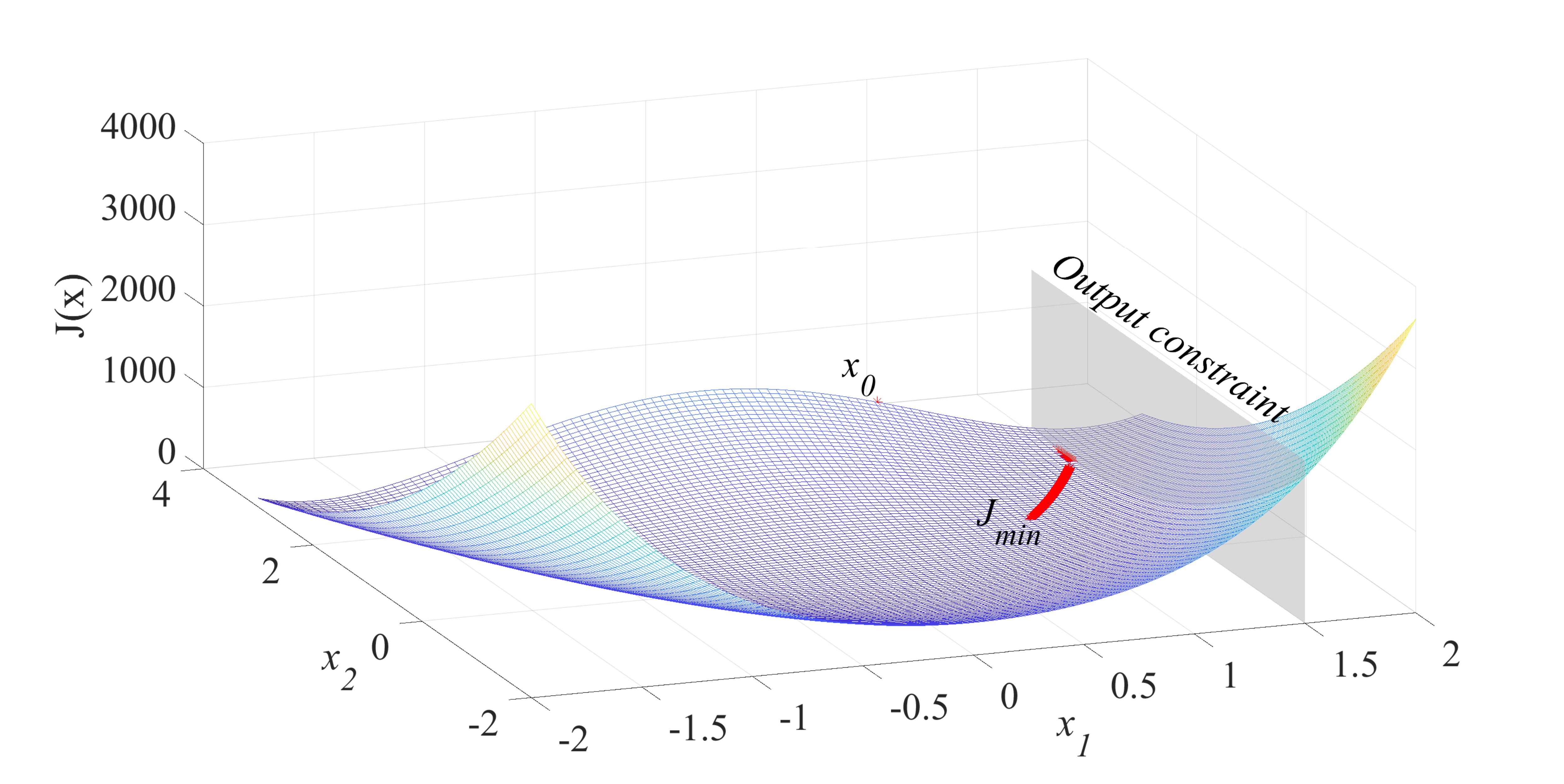}
	\caption{The  obtained state points via project gradient descent.}
	\label{fig:contourJ}
\end{figure}

Following Algorithm~\ref{alg:noesc}, also the simulation is divided into two steps: first, at each iteration, we solve the constrained optimization problem. Second, we design our controller for the finite-time state transition.  Given a initial guess $\boldsymbol{x}_{0} = [0.8 \, 3]^\top$, by making use of the project gradient descent in Algorithm~\ref{alg:grad_desc} with a fixed constant step $\alpha = 0.002$, a terminal condition $\varepsilon_0 = 10^{-2}$, and after $1524$ iterations, the minimum of the performance function $J_{\min}$ is obtained with acceptable tolerance and the resulting state sequence is shown in Fig.\ref{fig:contourJ} which satisfies the constraint $x_{1} \leq 1.5$.  The next step is to find the controller that can drive the system from the state $\boldsymbol{x}_k$ to $\boldsymbol{x}_{k+1}$. The first state transition, i.e., from $\boldsymbol{x}_{0} = [0.8~3]^\top$ to $\boldsymbol{x}_1 = [1.5~2.056]^\top$, is used to illustrate the proposed control design. The transition time interval $\Delta_k = 1$ is selected and the virtual output trajectory becomes for $0 \leq t \leq 1$
\begin{equation*}
	\zeta(t,p) = \zeta_{k} + \left(\zeta_{k+1} - \zeta_{k} - \gamma p \right) t + \gamma p t^2,
\end{equation*}
with boundary conditions 
\begin{equation*}
	\begin{aligned}
		\zeta_{0} = \frac{1}{4}\left[ \log \left(2.8\right) - \log \left( 1.2 \right) \right],~\zeta_{1} = \frac{1}{4}\left[ \log \left(3.5\right) - \log \left( 0.5 \right) \right],
	\end{aligned}
\end{equation*}
and the state of the internal dynamics is, for $0 \leq t \leq 1$,
\begin{equation*}
	\dot{\eta}^* = 2y^* - 2 \eta^{*3}, ~\eta^{*}_0 = 3,~ \eta^{*}_1 = 2.056.
\end{equation*}

\begin{figure}[h!]
	\centering
	\includegraphics[width=0.5\textwidth]{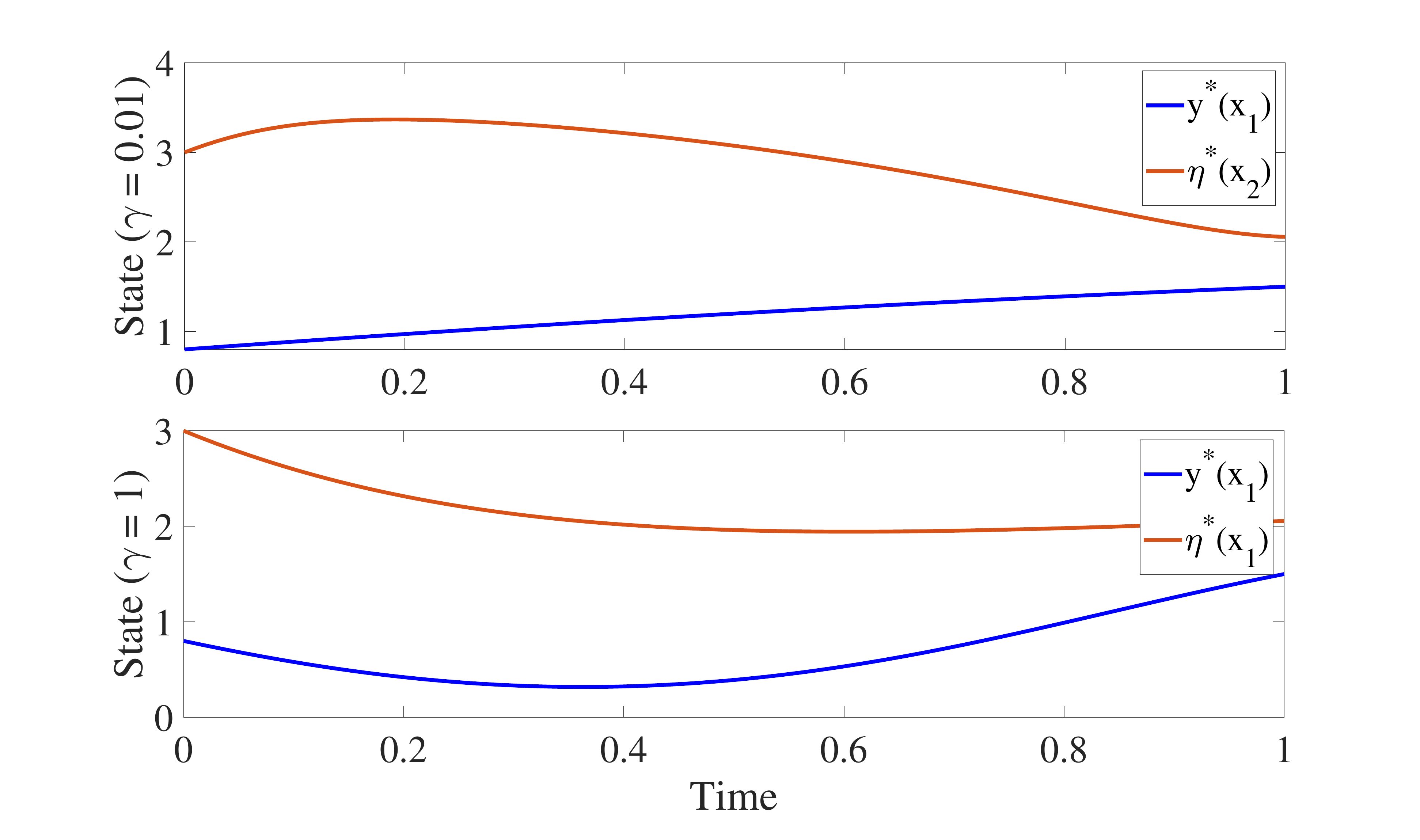}
	\caption{The resulting state trajectories with different selections of $\gamma$.}
	\label{fig:stateTraj1}
\end{figure}

\begin{figure}[h!]
	\centering
	\includegraphics[width=0.5\textwidth]{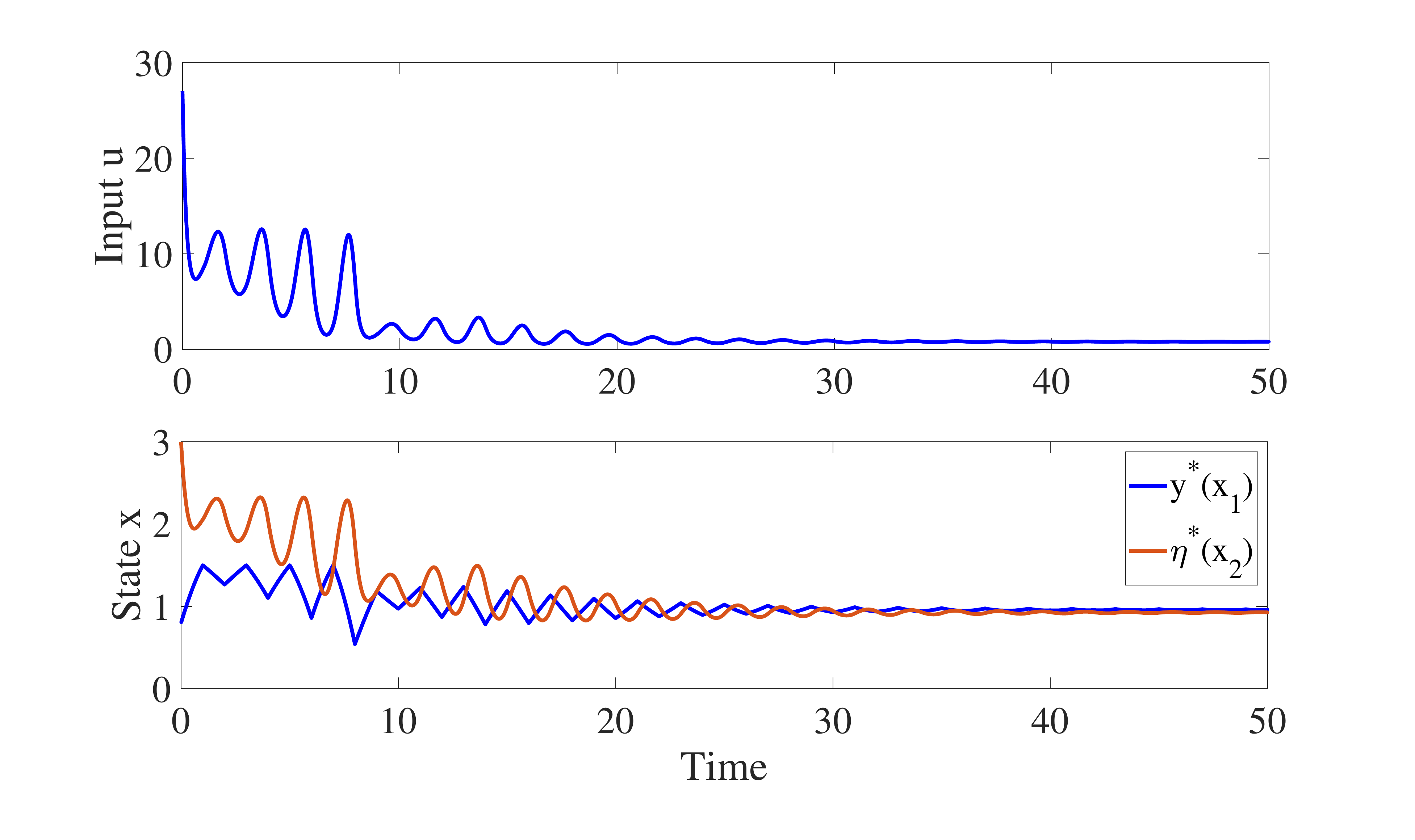}
	\caption{The seeking process with the resulting control input for $\rho = 1$ and $\gamma = 0.01$.}
	\label{fig:state_input}
\end{figure}

\begin{figure}[h!]
	\centering
	\includegraphics[width=0.5\textwidth]{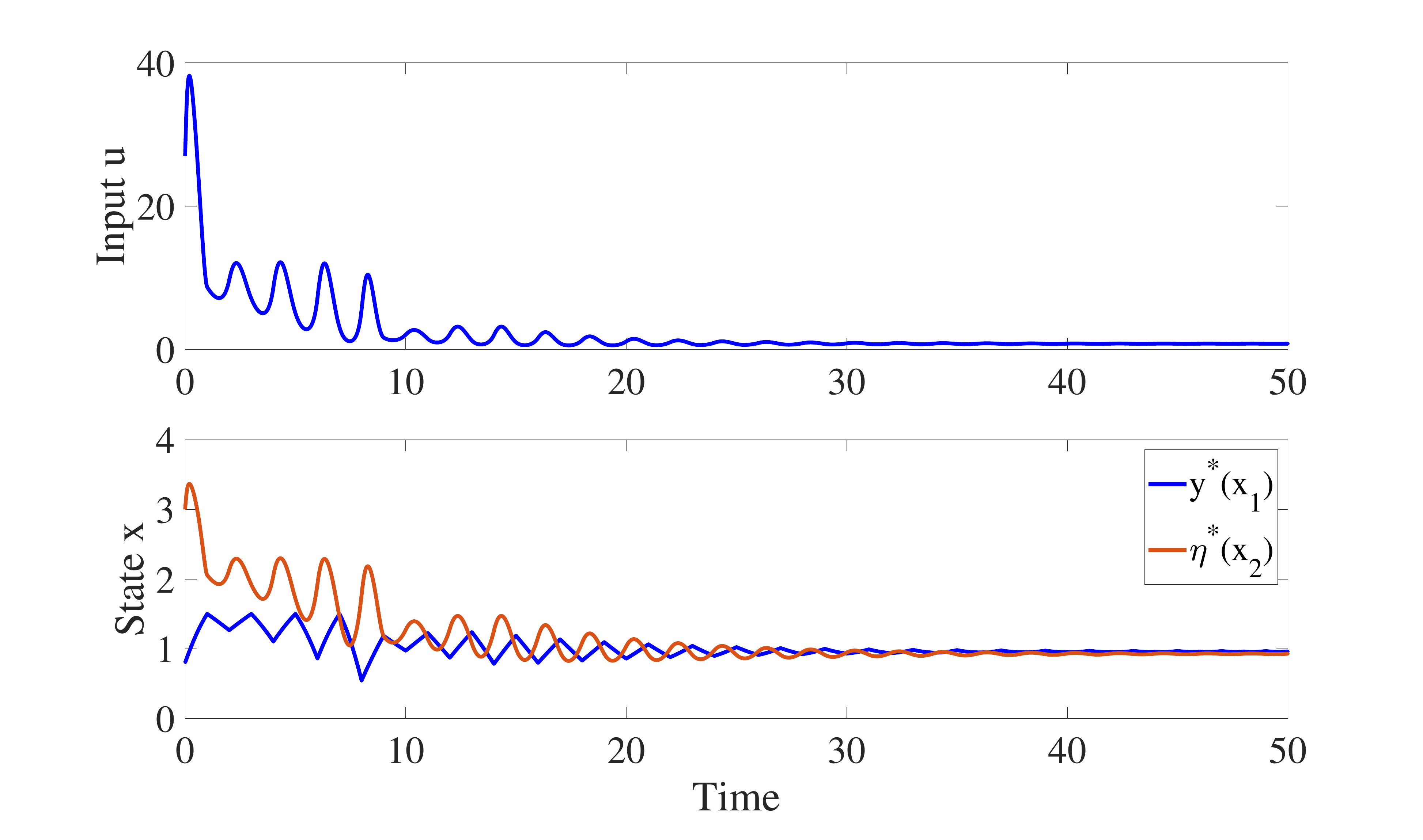}
	\caption{The seeking process with the resulting control input for $\rho = -1$ and $\gamma = 0.01$.}
	\label{fig:state_input2}
\end{figure}

\begin{figure}[h!]
	\centering
	\includegraphics[width=0.5\textwidth]{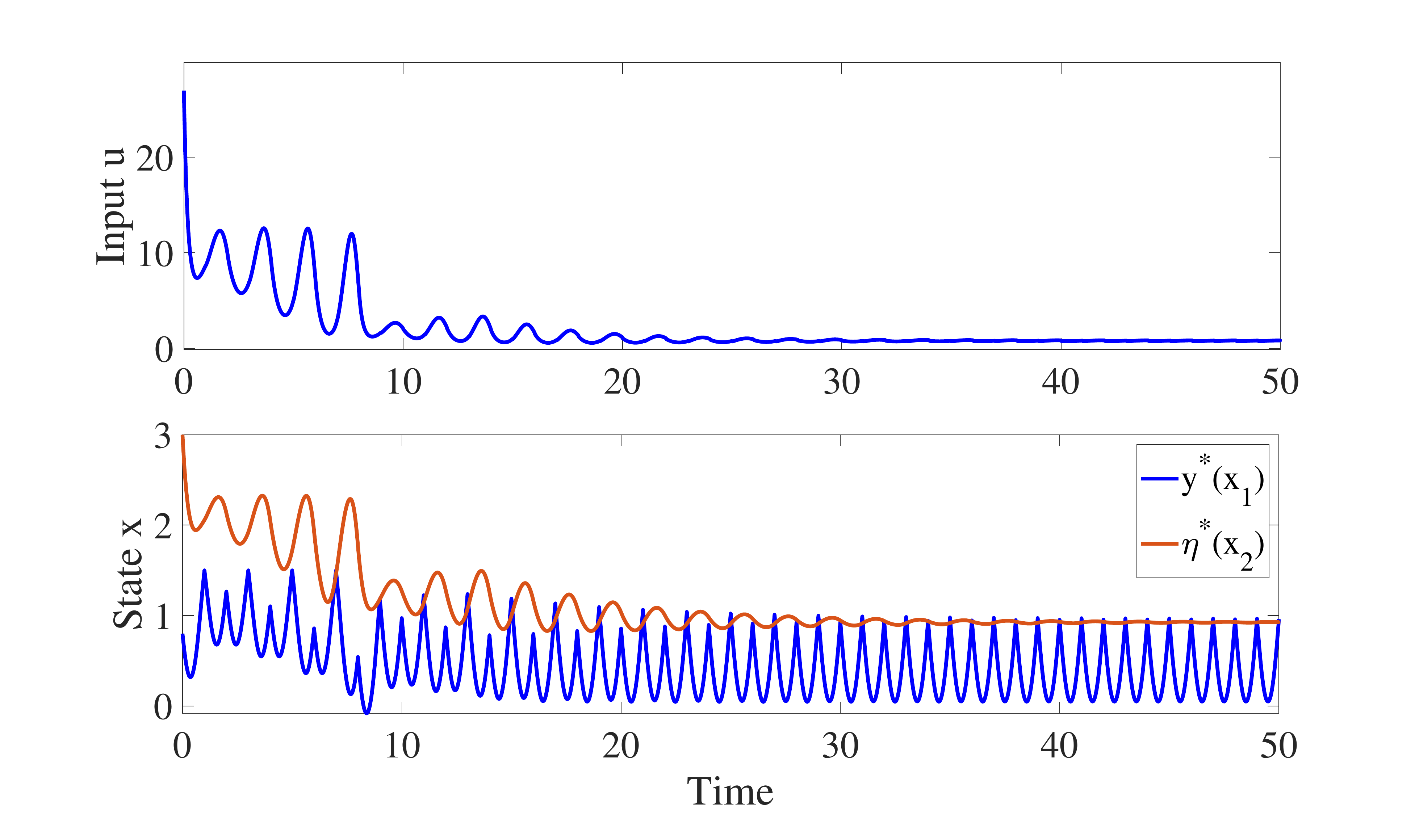}
	\caption{The seeking process with the resulting control input for $\rho = 1$ and $\gamma = 1$.}
	\label{fig:state_input3}
\end{figure}

To solve the boundary value problem, we adopt the MATLAB function \textit{bvp4c} asking to recast the boundary conditions of $\eta^*$ into the form of $\eta_{0}^{*}-3=0$ and $\eta_{1}^{*} - 2.056 = 0$. By selecting an initial guess for the parameter in \eqref{eq:zetaFunc} $p = 1$, an initial mesh of $100$ points, and the tuning parameter $\gamma = 0.01$, the resulting output and state trajectories of the state of the initial dynamics are shown in Fig.~\ref{fig:stateTraj1} (note that the first $50$ steps only have been explicitly shown in Fig.~\ref{fig:stateTraj1}--\ref{fig:state_input3}), which satisfy the boundary condition and thus the state transition from $\boldsymbol{x}_{0}$ to $\boldsymbol{x}_1$ is fulfilled by making use of the proposed control design. After solving the boundary condition problem for the whole seeking process $0\leq t \leq 1524$, the state trajectories and the corresponding input via the inversion-based approach are obtained as in Fig.~\ref{fig:state_input}, which shows that the states tend to the point $[1~1]^\top$, thus minimizing the performance function without violating the output constraint. 

For $\rho =-1$, the internal dynamics in \eqref{eq:IOsystExp} are unstable. With the proposed extremum seeking control design, the input and the resulting state trajectory are obtained as in Fig.~\ref{fig:state_input2}. Clearly, they show that the state variables converge to the stationary point $[1~1]^\top$ that minimizes the performance function, despite the fact that larger overshoot of $\eta$ in the first several steps is observed with respect to the case $\rho = 1$.  We have numerically supported the insight that the proposed \gls{NOESC} method can deal with nonlinear systems with unstable internal dynamics.  


\vspace{0.1cm}
\textit{B. Control performance with different selections of $\gamma$} 
\vspace{0.1cm}

To illustrate the impact of different selections of $\gamma$ on the state trajectories, we choose $\gamma = 1$ and by solving the two-point boundary problem, the resulting state trajectory from $x_{0}$ to $x_1$ is shown in Fig.~\ref{fig:stateTraj1}. Apparently, the output trajectory for $\gamma = 0.01$ is close to a straight line and has a smaller change of curvature compared to the case $\gamma = 1$, which implies better output performance with smaller value of $\gamma$ during the extremum seeking process. For the whole seeking process, the resulting state trajectory and input are shown in Fig.~\ref{fig:state_input3}, where a larger variation of the state $x_1$ appears due to larger value of $\gamma$. From our numerical experience, as a general rule of thumb, the smaller the value of $\gamma$, the better the state transition performance.

%%%%%%%%%%%%%%%%%%%%%%%%%%%%%%%%%%%%%%%%%%%%%%%%%%%%%%%%%%%%%%%%%%%%%%%%%
\section{Conclusions} 
In this paper, an extremum seeking control scheme for nonlinear systems subject to output constraints has been proposed. To find the extremum of an unknown performance function dependent on the state, constrained numerical optimization and finite-time state transition are combined. A novel NOESC framework that consists of the projected gradient descent and an inversion-based feedforward controller has been develop. The main merit of the proposed NOESC consists on the capabilities of deal with output constraints, and of finding the extremum of the performance function that directly depends on the states of the stable or unstable internal dynamics. Future work will focus on extending the results to the nonlinear system with state constraints. 

%%%%%%%%%%%%%%%%%%%%%%%%%%%%%%%%%%%%%%%%%%%%%%%%%%%%%%%%%%%%%%%%%%%%%%%%

\balance
\bibliographystyle{IEEEtran}
\bibliography{21_CDC_extremum_seeking}

\end{document}